\documentclass[a4paper,11pt,twoside]{article}

\usepackage{xypic}
\usepackage{latexsym}
\usepackage{amsmath}
\usepackage{amstext}
\usepackage{amsfonts}
\usepackage{amssymb}
\usepackage{amsbsy}
\usepackage{eucal}
\usepackage{enumerate}
\usepackage{theorem}

\xyoption{all}
\CompileMatrices

\newtheorem{theorem}{Theorem}

\newtheorem{lemma}{Lemma}
\newtheorem{proposition}[lemma]{Proposition}

\theorembodyfont{\rmfamily}
\newtheorem{remark}[lemma]{Remark}

\newcommand{\id}{\text{\rm id}}

\newcommand{\End}{\text{\rm End}}
\newcommand{\Aut}{\text{\rm Aut}}
\newcommand{\Gal}{\text{\rm Gal}}
\newcommand{\Mor}{\text{\rm Mor}}
\newcommand{\Hom}{\text{\rm Hom}}

\newcommand{\rank}{\text{\rm rank}}

\renewcommand{\char}{\text{\rm char}}

\newcommand{\F}{\mathbb{F}}

\newcommand{\FF}{\mathbb{F}}
\newcommand{\C}{\mathop{\cal C}}

\newcommand{\qed} {\hbox{} \nolinebreak \hfill $\;\Box$}

\begin{document}
\pagestyle{myheadings}
\markboth{\sc Bouw, Diem, Scholten}{\sc Ordinary elliptic curves}

\title{Ordinary elliptic curves of high rank over 
$\overline{\mathbb{F}}_{\!p}(x)$ with constant $j$-invariant}

\author{Irene I.\ Bouw
\and Claus Diem
\and Jasper Scholten}
\date{\today}

\maketitle

\begin{abstract}
We show that under the assumption of Artin's Primitive Root Conjecture,
for all primes $p$ there exist ordinary elliptic curves over
$\overline{\mathbb{F}}_{\!p}(x)$ with arbitrary high rank and constant
$j$-invariant. For odd primes $p$, this result follows from a theorem
which states that whenever $p$ is a generator of
$(\mathbb{Z}/\ell\mathbb{Z})^*/\langle -1 \rangle$ ($\ell$ an odd prime)
there exists a hyperelliptic curve over $\overline{\mathbb{F}}_{\!p}$
whose Jacobian is isogenous to a power of one ordinary elliptic curve.
\end{abstract}

\begin{center}\begin{minipage}{110mm}\small{\bf Key words~:}
      Elliptic Curves of High Rank, Jacobians. 
\end{minipage}
\end{center}

\begin{center}\begin{minipage}{110mm}\small{\bf 2000 MSC.~} 
Primary: 11G05; Secondary: 11G20, 14H40, 14H52. 
\end{minipage}
\end{center}

\section{Introduction}

Let $E$ be an elliptic curve over a field $L$. For 
various choices of $L$, it is known that $E(L)$ is a finitely generated 
group. This is the case if, for example,
\begin{itemize}
\item $L$ is a number field (by the Mordell-Weil Theorem, see
\cite{mordell}, \cite{weil}), or more generally
\item $L$ is finitely generated over its prime field (see 
\cite{neron}), or
\item $L$ is the function field of an algebraic variety
over a field $k$, and $E$ is not isogenous (over $L$) to an elliptic curve
which can be defined over $k$ (see \cite{langneron}).
\end{itemize}

One might ask how large the rank of $E(L)$ can get if one fixes $L$ and varies $E$. If $\char(L)=0$ then it is a well known open
problem whether this rank is
bounded or not in any of the above cases. But if $\char(L)$ is positive,
there are some results.
In the following table we list some cases for which it is known that the
rank can get arbitrary large.

{\small {\center \begin{tabular}{|c|c|c|c|} \hline
$L$ & $j$-\text{invariant} & \text{ord. / ss.} & authors\\
\hline \hline
$\overline{\F}_{\!2}(x)$ & constant & supersingular &
\begin{tabular}{c} well known,\\ see e.g.\ Elkies (1994), \cite{elkies}
\end{tabular}\\
\hline
\begin{tabular}{c} $\F_{\!p}(x)$ \\ $p$ odd \end{tabular}
& constant & supersingular & \begin{tabular}{c} Shafarevich and Tate\\ 
(1967), \cite{shaftate} \end{tabular} \\
\hline
$\overline{\F}_{\!2}(x)$ & constant &  ordinary & \begin{tabular}{c}
  follows from Gaudry, Hess and \\ Smart (2002) et al., \cite{GHS}, \cite{MQ}, \\
(see this paper) \end{tabular} \\
\hline
\begin{tabular}{c} $\overline{\F}_{\!p}(x)$ \\ $p$ odd \end{tabular}
& constant & ordinary &  \begin{tabular}{c}
B.D.S., assuming  Artin's \\ Primitive Root Conjecture \\ (see this paper) \end{tabular}
\begin{tabular}{c} \\ \end{tabular}
\\
\hline
$\F_{\!p}(x)$ & non-const. & & \begin{tabular}{c} \raisebox{0pt}[2.5ex]{Shioda (1986), \cite{shioda}, (for $\overline{\F}_{\!p}(x)$),} \\  Ulmer (2002), \cite{ulmer} \end{tabular} \\
\hline
\end{tabular}}}\\

Let $L$ be the function field of a (smooth, projective, geometrically irreducible) curve $C$ over
some field $k$ with $C(k) \neq \emptyset$. Let $E$ be an elliptic curve over $k$. 
It is well known
that there is a close relationship between $\rank(E(L)/E(k))$ and the
number of factors of $E$ in the Jacobian $J_C$ of $C$. Clearly,
$E(L) \simeq \Mor_k(C,E)$. Let some point $P \in C(k)$ be fixed. Then
$E(L)/E(k)$ is isomorphic to $\Mor_k((C,P),(E,\mathcal{O}))$, the
group of morphisms sending $P$ to the zero $\mathcal{O} \in
E(k)$. This group in turn is isomorphic to $\Hom_k(J_C,E)$. Let $J_C
\sim E^r \times A$ for some $r \in \mathbb{N}$ and some abelian
variety $A$ that does not have $E$ as a factor. Then the
$\mathbb{Q}$-vector space $\Hom_k^0(J_C,E)$ is isomorphic to
$\Hom_k^0(E^r,E) \simeq \Hom_k^0(E,E)^r$. So the rank of $E(L)/E(k)$
is equal to $r\cdot{}\rank(\End_k(E))$.

If $C$ is a hyperelliptic curve, and $k(x)$ is the rational quadratic subfield
of $L$, then one can consider the twist $E^{\rm twist}$ of $E_{k(x)}$ with respect to the extension $L|k(x)$. The action of the non-trivial element in $\text{Gal}(L|k(x))$ on $E(L) \! \otimes_{\mathbb{Z}} \! \mathbb{Q}$ induces a decomposition into eigenspaces
\[E(L) \otimes_{\mathbb{Z}} \mathbb{Q} = E(k(x)) \! \otimes_{\mathbb{Z}} \! \mathbb{Q} \, \oplus \, E^{\text{twist}}(k(x)) \! \otimes_{\mathbb{Z}} \! \mathbb{Q}.\]
Together with $E(k(x)) = E(k)$ this implies that $\rank(E^{\text{twist}}(k(x))) = $\linebreak $\rank(E(L)/E(k)) = r\cdot{}\rank(\End_k(E))$.
So one can construct high rank elliptic curves over $k(x)$ if one can
construct hyperelliptic curves over $k$ with a high factor $E^r$ in the
Jacobian. In \cite{shaftate} such curves are
given over prime fields $k$ of odd characteristic. These are supersingular and give
rise to the second line of the table. In \cite{elkies} this construction is
done over finite fields of characteristic 2, and the Mordell-Weil groups are studied in great detail.
The present paper deals with the case
of ordinary curves over finite fields.\\

In \cite{GHS} a new approach to attack the discrete-logarithm problem
in the group of rational points of an elliptic curve over a non-prime
finite field is given (see also \cite{He}, \cite{MQ}). The interest of
the authors of \cite{GHS} lies within the realm of cryptology but
their construction also gives rise to the following theorem which
implies the fourth line of the table (see Section \ref{sectionproof1}
for a proof). 

\begin{theorem}
\label{theorem-char-2}
For all $r \in \mathbb{N}$, there exists a hyperelliptic curve $H$
over $\mathbb{F}_{\!2^r}$ such that the Jacobian variety $J_H$ is
completely decomposable into ordinary elliptic curves and $J_H \sim
E^r \times A$ for some ordinary elliptic curve $E$ and a (ordinary,
completely decomposable) abelian variety $A$. If $r$ is a Mersenne
prime, there exists a hyperelliptic curve $H$ over
$\mathbb{F}_{\!2^r}$ of genus $r$ whose Jacobian variety is isogenous
to the power of one ordinary elliptic curve. 
\end{theorem}

In Section \ref{sectionproof2} of this paper, we prove the following theorem.

\begin{theorem}
\label{theorem-char-neq-2}
Let $p$ and $\ell$ be odd prime numbers such that $p$ generates \linebreak $(\mathbb{Z}/\ell\mathbb{Z})^*/\langle -1 \rangle$.
Then there exists a hyperelliptic curve $H$ over
$\overline{\mathbb{F}}_p$ 
of genus $\frac{\ell-1}{2}$ such that $J_H$ is isogenous to the power of one ordinary elliptic curve.
\end{theorem}

Recall that it is Artin's Primitive Root Conjecture that for a given non-square integer $a\neq -1$, there exist arbitrary large prime numbers $\ell$ with $\langle a\rangle=(\mathbb{Z}/{\!\ell}\mathbb{Z})^*$.
This conjecture
has not been proven for a single $a$. But it is known that there are at
most 2 prime values for $a$ for which Artin's Conjecture fails 
(see \cite{heathbrown}).
Also, it is proven that Artin's Conjecture follows from
the Generalized Riemann Hypothesis (see \cite{hooley}).

The fourth line of the table follows from Theorem
\ref{theorem-char-neq-2} and Artin's Conjecture for prime numbers $a$.

To the knowledge of the authors, it was not known before whether for
arbitrary large $r \in \mathbb{N}$ there exists some hyperelliptic curve over some field of characteristic
$\neq 2$ whose Jacobian variety is completely decomposable into
$r$ ordinary elliptic curves. The above Theorem \ref{theorem-char-neq-2}
also gives an affirmative answer to this question. Of
course, the question raised in \cite{ES} whether for all $r \in
\mathbb{N}$ there exist curves over $\mathbb{C}$ of genus $\geq r$
with completely decomposable Jacobian variety remains open.

\section{Proof of Theorem \ref{theorem-char-2}}
\label{sectionproof1}

We use the theory of function fields (in one variable) instead of the theory of curves. Let us fix the following notation: If $K$ is a perfect field and $L|K$ is a regular function field, we denote the Jacobian variety of the smooth, projective model of $L|K$ by $J_L$.

In the following, by a \emph{minimal subextension} of a field extension $\lambda|\kappa$ we mean some intermediate field $\mu$ of $\lambda|\kappa$ such that $\mu \supsetneq \kappa$ and $\mu|\kappa$ does not contain any non-trivial intermediate field.

We need the following lemma (see \cite{KR} and the proof of \cite[Theorem 2.1]{GS}).

\begin{lemma}
\label{elementary-abelian}
Let $K$ be a field, let $M|K(x)$ be a Galois extension with Galois group an elementary abelian $\ell$-group -- $\ell$ an arbitrary prime number -- such that $M|K$ is regular. Then $J_M \sim \prod_N J_N$ where $N$ runs over all minimal subextensions of $M|K(x)$. In particular, $g(M|K)$, the genus of $M|K$, is equal to $\sum_N g(N|K)$.
\end{lemma}

All the following extensions of $\mathbb{F}_{\!2}(x)$ should be regarded as embedded in a fixed algebraic closure $\overline{\mathbb{F}_{\!2}(x)}$. We use Artin-Schreier theory in the formulation of \cite[Theorem 8.3]{La}.

Fix some algebraic extension $K|\mathbb{F}_{\!2}$ and some $\alpha \in K \backslash \{ 0 \}$. Let $L|K(x)$ be the Artin-Schreier extension given by $y^2- y = x^{-1} + \alpha x$, i.e.\ $L$ corresponds by Artin-Schreier theory to the $\mathbb{F}_{\!2}$-vector subspace $\langle x^{-1} + \alpha x \rangle$ of $K(x)/\mathcal{P}(K(x))$, where $\mathcal{P} : K(x) \longrightarrow K(x), \; \xi \mapsto \xi^2 - \xi$ is the Artin-Schreier operator. Now $L|K$ is an ordinary elliptic function field -- the ordinarity follows for example from the Deuring-Shafarevich formula (see \cite[Corollary 1.8.]{Cr}) and the fact that $\overline{K}L|\overline{K}(x)$ has two ramified places --, and $J_L$ is an ordinary elliptic curve.

The action of the Galois group $\Gal(K|\mathbb{F}_{\!2}) \simeq
\Gal(K(x)|\mathbb{F}_{\!2}(x))$ on $K(x)$ gives rise to an action on
$K(x)/\mathcal{P}(K(x))$, and this action induces an action by
the group ring $\mathbb{F}_{\!2}[\Gal(K|\mathbb{F}_{\!2})]$. Let $U$ be the cyclic module generated by $x^{-1} + \alpha x$, and let $M|K(x)$ be the extension corresponding to $U$.

We claim that $M|K$ is regular. Note that the extension
$\overline{K}M|\overline{K}(x)$ is given by the image $\overline{U}$
of $U$ in $\overline{K}(x)/\mathcal{P}(\overline{K}(x))$, and $\overline{U}$ is isomorphic
to the image of $U$ in $K(x)/ \langle K \cup \mathcal{P}(K(x)) \rangle$. One sees easily that $U \longrightarrow \overline{U}$ is an isomorphism. It follows that $[M:K(x)] = [\overline{K}M : \overline{K}(x)]$, and $M|K$ is regular. 

The minimal subextensions $N$ of $M|K(x)$ all are either rational function fields or ordinary elliptic function fields. By Lemma \ref{elementary-abelian}, $J_M$ is an abelian variety which is completely decomposable into ordinary elliptic curves.

For some subextension $N$ of $M|K(x)$ and some $\sigma \in \Gal(K|\mathbb{F}_{\!2}) \simeq $\linebreak $\Gal(K(x)|\mathbb{F}_{\!2}(x))$, let $\sigma(N)$ be the image of $N$ in $M$ under some extension of $\sigma$ to $M$.

Let $V$ be the $\mathbb{F}_{\!2}$-vector subspace of $U$ which consists of the elements of the form $\beta x$ for some $\beta \in K$. Clearly, $[U:V] =2$. Let $R$ be the extension of $K(x)$ corresponding to $V$. Then by Lemma \ref{elementary-abelian}, the genus of $R$ is zero. Now, $[M:R] = [U:V] = 2$, thus $M$ is hyperelliptic.

Now let $r \in \mathbb{N}$. Let $\alpha \in \mathbb{F}_{\!2^r}$, not lying in any proper subfield, let $L$ and $M$ be defined as above with $K= \mathbb{F}_{\!2^r}$ and $\alpha$. Let $\sigma_{\mathbb{F}_{\!2^r}|\mathbb{F}_{\!2}} \in \text{Gal}(\mathbb{F}_{\!2^r}|\mathbb{F}_{\!2})$ be the Frobenius morphism. Then for $i = 0, \ldots, r-1$, the powers $\sigma_{\mathbb{F}_{\!2^r}|\mathbb{F}_{\!2}}^i(L)$ are pairwise distinct subfields of $M$. Now, all $J_{\sigma_{\mathbb{F}_{\!2^r}|\mathbb{F}_{\!2}}^i(L)}$ are isogenous to $J_L$ (via a power of the Frobenius homomorphism), and again by Lemma \ref{elementary-abelian}, $J_M \sim J_L^r \times A$ for some (ordinary, completely decomposable) abelian variety $A$ over $\mathbb{F}_{\!2^r}$.

It remains to prove the statement on the Mersenne primes.

Let $r \in \mathbb{N}$ be an odd prime. Let $\beta$ be a generator of
the $\mathbb{F}_{\!2}[\text{Gal}(\mathbb{F}_{\!2^r}|\mathbb{F}_{\!2})]$-module
$\mathbb{F}_{\!2^r}$ (i.e.\ $\beta$,
$\sigma_{\mathbb{F}_{\!2^r}|\mathbb{F}_{\!2}}(\beta), \ldots, 
\sigma_{\mathbb{F}_{\!2^r}|\mathbb{F}_{\!2}}^{r-1}(\beta)$ form a normal
basis of $\mathbb{F}_{\!2^r}|\mathbb{F}_{\!2}$). 

Let $\varphi_2(r)$ be the (multiplicative) order of 2 modulo $r$.

Recall that we have canonical isomorphisms $\mathbb{F}_{\!2}[\text{Gal}(\mathbb{F}_{\!2^r}|\mathbb{F}_{\!2})] \simeq \mathbb{F}_{\!2}[\mathbb{Z}/r\mathbb{Z}] \simeq \mathbb{F}_{\!2}[x]/(x^r-1)$ of rings, and we have a decomposition into irreducible factors $x^r-1 = (x-1) p_1 \cdots p_{\frac{r-1}{\varphi_2(r)}}$ where the $p_i$ are pairwise distinct polynomials of degree $\varphi_2(r)$.

Let $\alpha:=(((x-1)p_2 \cdots
p_{\frac{r-1}{\varphi_2(r)}})(\sigma_{\mathbb{F}_{\!2^r}|\mathbb{F}_{\!2}}))(\beta)$.
Then $\alpha \notin \mathbb{F}_{\!2}$,
$(p_1(\sigma_{\mathbb{F}_{\!2^r}|\mathbb{F}_{\!2}}))(\alpha) = 0$. Let $M$ be
defined as above. Now the assignment $f \mapsto
\mathbb{F}_{\!2^r}(x)[\mathcal{P}^{-1}(x^{-1} +
f(\sigma_{\mathbb{F}_{\!2^r}|\mathbb{F}_{\!2}})(\alpha) \, x)]$ induces a
bijection between the polynomials $\neq 0$ of degree $< \deg(p_1)$ and
the minimal subextensions $N$ of $M|\mathbb{F}_{\!2^r}(x)$ with genus
1. There are $2^{\varphi_2(r)}-1$ such polynomials, and thus the genus
of $M$ is $2^{\varphi_2(r)}-1$.

Now let $r$ be a Mersenne prime, i.e.\ $r$ is a prime of the form
$2^\ell-1$. Then $\varphi_2(r) =\ell$ and there are $2^{\varphi_2(r)}-1 = r$
minimal subextensions $N$ of $M|\mathbb{F}_{\!2^r}(x)$ of genus 1. These subextensions are equal to the
$r$ distinct subextensions
$\sigma_{\mathbb{F}_{\!2^r}|\mathbb{F}_{\!2}}^i(L)|\mathbb{F}_{\!2^r}(x)$ of
genus 1. Thus $J_M \sim J_{L}^r$.
\qed

\section{Proof of Theorem \ref{theorem-char-neq-2}}
\label{sectionproof2}

The idea of the proof of Theorem \ref{theorem-char-neq-2} is to
consider curves $C$ over a finite field $k$ such that after some base
extension $K|k$, $J_{C_K}$ has an endomorphism not defined over any
proper subextension of $K|k$. If additionally $J_C$ is ordinary,
this endomorphism induces a decomposition of $J_{C_K}$ as is made
precise in the next subsection. 

We then apply this general result to hyperelliptic curves in certain algebraic families. These families have already been studied in characteristic 0 in \cite{ttv}. We use techniques similar to those of \cite{Bo} to show that they are generically ordinary.

\subsection{Operation on abelian varieties over finite fields}
In this subsection, we deal with the following situation:

Let $K|k$ be an extension of finite fields inside the fixed algebraic closure $\overline{k}$ of $k$. Let $\sigma_{K|k} \in \text{Gal}(K|k)$ be the Frobenius morphism, and let $\ell \neq \text{char}(k)$ be a prime. Let $A$ be an abelian variety over $k$.

Assume furthermore that we are given an $\tau \in \End_K^0(A_K)$ such that
\begin{enumerate}
\item
the action of $\Gal(K|k)$ on $\End^0_K(A_K)$ restricts to $\mathbb{Q}[\tau] \leq \End^0_K(A_K)$,
\item
$\tau$ is not defined over any intermediate field $\mu$ of $K|k$ with $\mu \subsetneq K$,
\item
$\mathbb{Q}[\tau]$ is a field.
\end{enumerate}

\begin{proposition}
\label{char-poly-structure}
Under the above assumptions, the characteristic polynomial of the Frobenius endomorphism of $A$ has
the form $f(T^{[K:k]})$ for some polynomial $f(T) \in \mathbb{Z}[T]$
of degree $\frac{2\dim(A)}{[K:k]}$.
\end{proposition}

As a special case of this proposition, we obtain.
\begin{proposition}
\label{Weil-res}
If additionally to the above assumptions $A$ is ordinary then $A$ is isogenous to the Weil restriction with
respect to $K|k$ of an ordinary abelian variety $B$ over $K$ with
$\dim(B) = \frac{\dim(A)}{[K:k]}$. In particular, $A_K \sim
B^{{\dim(A)}/{[K:k]}}$.
\end{proposition}
\emph{Proof of Proposition \ref{Weil-res} assuming Proposition \ref{char-poly-structure}.}  Let $\chi_A$ be the characteristic polynomial of the Frobenius endomorphism of $A$. By Proposition
\ref{char-poly-structure}, $\chi_A = f(T^{[K:k]})$ for some polynomial
$f \in \mathbb{Z}[T]$ of degree $\frac{2\dim(A)}{[K:k]}$. This implies
that $\chi_{A_K} = f(T)^{[K:k]}$. 

There exists an (ordinary) abelian variety $B$ over $K$ such that $\chi_B = f$. (For every irreducible factor $f_i$ of $f$, there exists some $K$-simple abelian subvariety $B_i$ of $A_K$ such that $\chi_{B_i}$ is a power of $f_i$. As $A_K$ is ordinary by assumption, so is $B_i$. This implies that $\chi_{B_i}$ is irreducible, and consequently that $\chi_{B_i} = f_i$.)

The Weil restriction of $B$
with respect to $K|k$ has characteristic polynomial $\chi_B(T^{[K:k]})
= \chi_A$ (see \cite[\S1 (a)]{Mi}). This implies that $A \sim
\text{Res}_k^K(B)$ (see \cite[Appendix 1, Theorem 2]{Mu}).  \qed
\\

\emph{Proof of Proposition \ref{char-poly-structure}.}  By assumption,
the action of the Galois group $\Gal(K|k)$ on $\mathbb{Q}[\tau]$ gives an injective
homomorphism $\Gal(K|k) \longrightarrow \Aut(\mathbb{Q}[\tau])$. Fix
some polynomial $p(T) \in \mathbb{Q}[T]$ such that $\sigma_{K|k}(\tau)
= p(\tau)$. For $i \in \mathbb{N}_0$, define $p_i$ by $p_0 := T, \; 
p_{i+1} := p_i(p(T))$. Then $\sigma^{\,i}_{K|k}(\tau) =
p_i(\tau)$. This implies that the elements $p_i(\tau)$ for $i=0,
\ldots, [K:k]-1$ are pairwise distinct and $p_{[K:k]}(\tau)=\tau$.

Let $V_\ell(A) := T_\ell(A) \otimes_{\mathbb{Z}_\ell} \mathbb{Q}_\ell$,
$\overline{V}_{\!\!\ell}(A) := V_\ell(A) \otimes_{\mathbb{Q}_\ell}
\overline{\mathbb{Q}}_\ell$.

We will show that the characteristic polynomial of the Frobenius endomorphism in
its operation on $T_\ell(A)$ (or -- what amounts to the same -- on
$\overline{V}_{\!\!\ell}(A)$) has the form $f(T^{[K:k]})$ for some polynomial
$f(T) \in \overline{\mathbb{Q}}_\ell[T]$ of degree
$\frac{2\dim(A)}{[K:k]}$. As $f(T^{[K:k]}) \in \mathbb{Z}[T]$, the same holds for $f(T)$.

As by assumption $\mathbb{Q}[\tau]$ is a field, the operation of $\tau$ on $\overline{V}_{\! \ell}$ is diagonalizable.
For some eigenvalue $\lambda$ of $\tau$ in its operation on
$\overline{V}_{\!\!\ell}(A)$, let $\overline{V}_{\!\!\ell}^\lambda$ be the corresponding eigenspace.

Let $\pi_k$ be the Frobenius endomorphism of $A$ over $k$. Then for $P
\in A(\overline{k})$, we have $\pi_k(P) = \sigma_k^{-1}(P)$,
where $\sigma_k \in \Gal(\overline{k}|k)$ is the Frobenius morphism. This
implies $\alpha \pi_k = \pi_k \, \sigma_{K|k}(\alpha)$ for all $\alpha
\in \End_K^0(A_K)$, thus
\begin{equation}
\label{pi-tau}
\tau \pi_k^{\, i} = \pi_k^{\, i} \,
\sigma_{K|k}^{\,i}(\tau) = \pi_k^i \, p_i(\tau) \; \text{ for } i \in
\mathbb{N}_0.
\end{equation}
Fix some eigenvalue $\lambda$ and some $i \in \mathbb{N}$. Then by equation (\ref{pi-tau}), $\pi_k^i(\overline{V}_{\!\!\ell}^\lambda) \leq
\overline{V}_{\!\!\ell}^{p_i(\lambda)}$. (In particular, $p_i(\lambda)$ is
an eigenvalue of $\tau$.) Since $\overline{V}_{\!\!\ell}(A)$ is the direct sum of the eigenspaces for
$\tau$ and $\pi$ is bijective, we have
\[\pi_k^i(\overline{V}_{\!\!\ell}^\lambda) = \overline{V}_{\!\!\ell}^{p_i(\lambda)}.\]

The equation $p_{[K:k]}(\tau) = \tau$ implies that $p_{[K:k]}(\lambda)
= \lambda$. We claim that the eigenvalues $\lambda = p_0(\lambda), \; p(\lambda) = p_1(\lambda), \; \ldots, \; p_{[K:k]-1}(\lambda)$ are pairwise distinct.

To prove this, note that $\lambda$ is a root of $\chi_\tau = m$, thus $\mathbb{Q}[\lambda] \simeq \mathbb{Q}[T]/(m(T)) \simeq \mathbb{Q}[\tau]$. The claim on the eigenvalues follows from the fact that the $p_i(\tau)$ are pairwise distinct for $i=0, \ldots, [K:k]-1$. 

We have a direct sum $\bigoplus_{i=0}^{[K:k]-1}(\overline{V}_{\!\!\ell}^{p_i(\lambda)}) \leq \overline{V}_{\!\!\ell}(A)$ which we denote by $\overline{V}_{\!\!\ell}(\lambda)$.
The operation of $\pi_k$ on $\overline{V}_{\!\!\ell}(A)$ restricts to $\overline{V}_{\!\!\ell}(\lambda)$, and on this $[K:k] \cdot \dim(\overline{V}_{\!\!\ell}^{\lambda})$-dimensional space, $\pi_k$ can be described by a block matrix of the form
\[ \left(\begin{array}{ccccc}
O &     &        & M_\lambda \\
I & O          &        &   \\
  &     \ddots & \ddots &   \\
  &            & I      & O
\end{array} \right),\]
where each of the blocks $O$, $I$, $M_\lambda$ has dimension $\dim(\overline{V}_{\!\!\ell}^{\lambda}$).

One sees that on $\overline{V}_{\!\!\ell}(\lambda)$, the characteristic polynomial of the Frobenius endomorphism has the desired form. The result follows from the fact that $\overline{V}_{\!\!\ell}(A) = \bigoplus_\lambda \overline{V}_{\!\!\ell}(\lambda)$, where $\lambda$ runs over a certain subset of the set of eigenvalues of $\tau$.
\qed

\subsection{Some families of hyperelliptic curves}

In this subsection, we want to study the $p$-rank of curves in certain
families of hyperelliptic curves.

Let $p$ be an odd prime. For a field $k$ of characteristic $p$, a $t
\in k \backslash \{\pm 2\}$ and an odd $\ell$ prime to $p$, let
$C^\ell_t$ (or $C_t$ if $\ell$ is fixed) be the hyperelliptic curve
over $k$ given the affine equation
\[y^2 = x(x^{2\ell} + tx^\ell + 1). \]
The goal of this subsection is to prove the following proposition.

\begin{proposition}
\label{ordinaryprop}There exists an open subset $U\subset \mathbb{A}_{\mathbb{F}_{\!p}}^1 \backslash \{\pm 2\}$ such that
\begin{enumerate}[(a)]
\item 
for every $\ell$ as above, every field $k$ of characteristic $p$ and every $t \in U(k)$, the curve $C^{\ell}_t$ is ordinary, 
\item
if $i \in \mathbb{N}, i > 1$, then $U(\FF_{\!p^i})$ is nonempty.
\end{enumerate}
\end{proposition}

Fix some $\ell$, some perfect field $k$ containing the $\ell$th roots of unity and
$t\in k \backslash \{\pm 2\}$. Choose a primitive $2\ell$th root of unity
$\zeta_{2\ell} \in k$ and define an automorphism $\tau_{2\ell}$ of
$C^\ell_t$ by $(x,y) \mapsto (\zeta_{2\ell}^2x,
\zeta_{2\ell} y)$.

Note that the genus of $C^\ell_t$ is $\ell$. The holomorphic
differentials $\omega_i$ defined by
\[
\omega_i= x^{i-1} \frac{{\rm d}x}{y}, \qquad i=1, \ldots,\ell
\]
form a basis of $H^0(C^\ell_t, \Omega)$ (see \cite{Yui}). Moreover,
$\tau_{2\ell}\, \omega_i=\zeta_{2\ell}^{2i-1}\omega_i$. Therefore
$\omega_i$ is an eigenvector of $\tau_{2\ell}$ with eigenvalue
$\zeta_{2\ell}^{2i-1}$. 

The Cartier operator $\C:H^0(C^\ell_t, \Omega)\to
H^0(C^\ell_t, \Omega)$ is defined as the dual with respect to Serre
duality of the absolute Frobenius $F:H^1(C^\ell_t, {\cal O})\to
H^1(C^\ell_t, {\cal O})$. It is $\mathbb{F}_{\!p}$-linear and satisfies $\mathcal{C}\alpha^p \omega = \alpha \, \mathcal{C} \omega \; \; (\alpha \in k, \, \omega \in H^0(C^\ell_t, \Omega))$. It is a bijection if and only if $C^\ell_t$
is ordinary. We want to describe the matrix of $\C$ with respect to the above basis of $H^0(C^\ell_t, \Omega)$. In order to do so, we need some more notation.

For $i\in\{1,\ldots, \ell\}$,
define $j(i)\in\{1,\ldots, \ell\}$ and $\alpha(i)\in\{0, \ldots,
p-1\}$ by
\[
2j(i)-1\equiv \frac{2i-1}{p} \pmod{2\ell}, \qquad
\alpha(i)=\left[\frac{p(2j(i)-1)}{2\ell}\right].
\] 
Here $[\cdot]$ denotes the integral part, as usual.

Let $f:=(x^2+tx+1)^{(p-1)/2} \in \mathbb{F}_{\!p}[t,x]$ and write $f= \sum_{n=0}^{p-1} c_n x^n$ with $c_n \in \mathbb{F}_{\!p}[t]$. Note that
\begin{eqnarray*}
\label{binomial}
c_n&=&\sum_{2n_1+n_2=n}\binom{(p-1)/2}{n_1}\binom{(p-1)/2 - n_1}{n_2}t^{n_2}.
\end{eqnarray*}
  For later use we
remark that if $n \leq \frac{p-1}{2}$, then $\deg(c_n) = n$ (because
$\binom{(p-1)/2}{n} \neq \nolinebreak 0$).

Now let $k:=\overline{\mathbb{F}_{\!p}(t)}$ and let $C_t^\ell$ be defined as above.
\begin{lemma}\label{coeflem}
For every $i\in\{1, \ldots, \ell\}$, we have
\[
\C \omega_i= c_{\alpha(i)}^{1/p}\, \omega_{j(i)}.
\]
\end{lemma}
\emph{Proof.}
As the automorphism $\tau_{2\ell}$ commutes with the absolute Frobenius $F$ in their operation on $H^1(C,\mathcal{O})$, the operation of $\tau_{2\ell}$ on $H^1(C,\Omega)$ commutes with $\mathcal{C}$. This implies that $\C\omega_i$ is an
eigenvector of $\tau_{2\ell}$ with eigenvalue
$\zeta_{2\ell}^{(2i-1)/p}$. In particular,
$\C\omega_i=\gamma_i^{1/p}\omega_{j(i)}$, for some $\gamma_i\in
k$. We want to show that $\gamma_i=c_{\alpha(i)}$. 

The Cartier operator extends to an $\mathbb{F}_{\!p}$-linear operator $\mathcal{C}$ on the meromorphic differentials which satisfies $\mathcal{C}h \omega = h^p \, \mathcal{C} \omega\; \;(h \in k(C^\ell_t), \, \omega \in \Omega(k(C^\ell_t))$. It is well known that
$\C \frac{{\rm d} x}{x}=\frac{{\rm d}x}{x}$ and $\C x^i{\rm d}x=0 \mbox{ if } p\nmid (i-1)$
(see for example \nolinebreak \cite{Yui}).  

We have
\[
\omega_i=x^{i-1}\frac{{\rm d} x}{y}=\frac{x^{p j(i) }}{y^p} x^{(p-1)/2+i-j(i)p}
  f(x^\ell)\frac{{\rm d} x}{x}.
\]
Define $g=x^{(p-1)/2+i-p j(i)}f(x^\ell)$ and write $g=\sum_m g_m x^m.$
Then
\[
\C\omega_i=\frac{x^{j(i)-1}}{y}\left[\sum_m g_{pm}^{1/p}x^m\right]{\rm d}x.
\]
We want to find all $m$ such that $g_{pm}\neq 0$. The definition of
$g$ implies that $g_{pm}=c_n$, where 
\[
pm=\frac{p-1}{2} +i-p j(i)+n\ell.
\]
Recall that the degree of $f$ is $p-1$.
Therefore, we need to find all $n$ such that $0\leq n\leq p-1$ and
\begin{equation}
\label{coeff-neq-0}
p-1+2i-2pj(i)+2n\ell \equiv 0 \pmod{p}.
\end{equation}
Because of the equality
\[
p(2j(i)-1)=2\ell\langle\frac{p(2j(i)-1)}{2\ell}\rangle+2\ell\left[\frac{p(2j(i)-1)}{2\ell}\right]=(2i-1)+\alpha(i)2\ell,
\]
(\ref{coeff-neq-0}) equivalent to $2n\ell \equiv 2\ell\alpha(i) \pmod{p}$. The only such $n$ is $n=\alpha(i)$. This proves the lemma.  \qed

\bigskip\noindent  
\emph{Proof of Proposition \ref{ordinaryprop}.}  Let $\ell$,
$k=\overline{\mathbb{F}_{\!p}(t)}$ and $C^\ell_t$ be as above.  Let $A^{(\ell)}$ be the
matrix obtained by raising all coefficients of the matrix of the
Cartier operator to the $p$th power.  Lemma \ref{coeflem} shows that
$A^{(\ell)}$ is the product of a permutation matrix
and the diagonal matrix $(c_{\alpha(i)} \delta_{i,j})_{i,j}$, where $\delta_{i,j}$ is the Kronecker delta. (Note that the $\alpha(i)$ depend on $\ell$.) Define
\[
\Phi:=\prod_{n=0}^{(p-1)/2} c_{n}.
\]
Since $c_n=c_{p-1-n}$, the determinant of $A^{(\ell)}$ divides a sufficiently
large power of $\Phi$. 

Now let $k$ be an arbitrary perfect field of characteristic $p$, and choose some $t_0 \in k \backslash \{ \pm 2 \}$. Analogous to above, let $A_{t_0}^{(\ell)}$ be the matrix obtained by raising all coefficients of the matrix of the
Cartier operator of $C_{t_0}^\ell$ to the $p$th power. Then $A_{t_0}^{(\ell)}$ is the specialization of $A_t^{(\ell)}$ induced by the homomorphism $\mathbb{F}_{\!p}[t] \longrightarrow k, \; t \mapsto t_0$.

This implies that the curve $C^\ell_{t_0}$ is
ordinary if $\Phi(t_0)\neq 0$.  Now define $U:=\mathbb{A}_{\mathbb{F}_{\!p}}^1 \backslash (\{\pm 2\}\cup\{t\, |\, \Phi(t)=0\})$. Obviously $U$ does not depend on $\ell$.

We have already seen that $\deg(c_n)=n$ for $n \leq \frac{p-1}{2}$. Therefore
\[
\deg(\Phi)=\sum_{n=0}^{(p-1)/2} n = \frac{p^2-1}{8} < p^2-2.
\]
This proves (b).
\qed

\subsection{Completely decomposable Jacobians}
Fix some distinct odd prime numbers $p$ and $\ell$. For a field $p$ of characteristic $p$ and a $t \in k$, let $E_t$ be the elliptic
curve given by the affine equation
\[y^2 = x(x^2+tx+1).\]
We have a cover  $\pi: C_t \longrightarrow
E_t, \;\; (x,y) \mapsto (x^\ell,y \, x^{(\ell-1)/2})$.

Let $k:=\mathbb{F}_q$, where $q$ is some power of $p$, and
choose some $t\in k \backslash \{\pm 2\}$.  Let $K:=\mathbb{F}_{\!q}[\zeta_\ell]$.

Let $A_t$ be the reduced identity component of the kernel of $\pi_* : J_{C_t} \longrightarrow J_{E_t}$ -- this is an $(\ell -1)$-dimensional abelian variety. It is equal to the complement under the canonical principal polarization of $J_{C_t}$ of $\pi^*(J_{E_t})$.

Let $\tau_\ell := \tau_{2\ell}^2$. We have $\pi^*(J_{E_t}) = \pi^*
\pi_*(J_{C_t}) = (1 + \tau^*_\ell + \cdots + {\tau^*_\ell}^{\ell-1})(J_{C_t})$, and $A_t = (1 - \frac{1 + \tau^*_\ell  + \cdots + {\tau_\ell^*}^{\ell-1}}{\ell})(J_{C_t})$. (Note that $1 + \tau^*_\ell + \cdots + {\tau_\ell^*}^{\ell-1}$ is invariant under the Galois action and thus lies in $\End_k^0(J_{C_t})$.)

This implies:
\begin{lemma}
The automorphism $\tau_{\ell}^*$ restricts to a $K$-automorphism of $(A_t)_K$, and $\mathbb{Q}[\tau_{\ell}] \leq \End^0_K((A_t)_K)$ is a field (isomorphic to $\mathbb{Q}[\zeta_{\ell}]$, where $\sigma_{K|k} \in \Gal(K|k)$ operates by $\zeta_{\ell} \mapsto \zeta_{\ell}^q$).
\end{lemma}

Now let $i \in \mathbb{N}, i>1$ and assume that $p^i$ is a generator modulo $\ell$. By Proposition
\ref{ordinaryprop}, there exists some $t \in \mathbb{F}_{\!p^i}
\backslash \{\pm 2 \}$ such that $C_t$ and thus $J_{C_t}$ is ordinary.

Again let $k:= \mathbb{F}_{\!p^i}$, $K:= k[\zeta_\ell]$. Then $\tau^*_\ell|_{(A_t)_K}$ is not defined over any subfield $\mu$ of $K|k$ with $\mu \subsetneq K$ and $[K:k]=\ell-1=\dim(A_t)$. We can thus apply Proposition \ref{Weil-res} to $A_t$, $K|k$ and $\tau_\ell^*$.

We conclude that $A_t$ is the Weil restriction (with respect to $K|k$)
of an ordinary elliptic curve over $K$. It follows that $J_{C_t} \sim
E_t \times \text{Res}^K_k(\widetilde{E}_t)$ for some elliptic curve
$\widetilde{E}_t$ over $K$. This implies that $J_{{(C_t)}_{K}} \sim
  (E_t)_K \times (\widetilde{E}_t)_K^{\ell-1}$.

We have proven:
\begin{proposition}
Let $p$ and $\ell$ be odd prime numbers and $i \in \mathbb{N}, i > 1$,
such that $p^i$ is a (multiplicative) generator modulo $\ell$. Then
there exists a hyperelliptic curve over $\mathbb{F}_{\!p^{i}}$ of genus
$\ell$ whose Jacobian variety becomes over $\mathbb{F}_{\!p^{i(\ell-1)}}$
isogenous to the product of one ordinary elliptic curve and the $(l-1)$th power
of one ordinary elliptic curve.
\end{proposition}

This proposition already implies the fourth line of the table in the introduction. In order to prove Theorem \ref{theorem-char-neq-2}, let us study the
hyperelliptic curves $C_t$ ($k$ arbitrary, $t \in k \backslash \{\pm 2\}$) in more detail.

In addition to the automorphism $\tau_{2\ell}$,
$C_t$ has the automorphism $\gamma :  (x,y) \mapsto (\frac{1}{x},\frac{y}{x^{\ell+1}})$ of order 2.
Let $D_t$ be the quotient of $C_t$ by this automorphism,
$c : C_t \longrightarrow D_t$ the covering morphism. $D_t$ is given by
the equation
\[y^2 = D_\ell(x,1) + t, \]
where $D_\ell(x,a) := (\frac{x+\sqrt{x^2-4a}}{2})^\ell + (\frac{x-\sqrt{x^2-4a}}{2})^\ell \in k[x]$ is the $\ell$-the Dickson polynomial for $a \in k^*$ (c.f.\ \cite{dickson}). With this equation, $c: C_t \longrightarrow D_t$ is given by $(x,y) \mapsto (x + x^{-1}, \frac{y}{x^{(\ell+1)/2}})$. This follows from the equation
\[D_\ell(x + \frac{a}{x},a) = x^\ell + (\frac{a}{x})^\ell.\]
We see in particular that $D_t$ has genus $\frac{\ell-1}{2}$. Note also
that if $C_t$ is ordinary so is $D_t$. Thus in particular, if $i > 1$
there exists some $t \in \F_{\!p^i}$ such that $D_t$ is ordinary.

The covering morphism $c: C_t \longrightarrow D_t$ induces
canonical homomorphisms $c^* : J_{D_t} \longrightarrow
J_{C_t}$ and $c_* : J_{C_t} \longrightarrow
J_{D_t}$. The following argument shows that the kernel of $c_* :
J_{C_t} \longrightarrow J_{D_t}$ contains $\pi^*(E_t)$, and the image of $c^*
: J_{D_t} \longrightarrow J_{C_t}$ is contained in $\ker(\pi_*) = A_t$.

We have the identity
$\gamma \tau_\ell = \tau^{-1}_\ell \gamma$ in $\Aut(C_t)$. This identity implies $(\id + \tau^* +
\cdots + {\tau^*_\ell}^{\ell-1}) \gamma^* = \gamma^* (\id + \tau_\ell^* +
\cdots + {\tau_\ell^*}^{\ell-1})$ on $J_{C_t}$. This in turn implies that both $A_t = \ker(\id + \tau_\ell^* +
\cdots + {\tau_\ell^*}^{\ell-1})$ and $\pi^*(E_t) = (\id + \tau_\ell^* +
\cdots + {\tau_\ell^*}^{\ell-1})(J_{C_t})$ are invariant under $\gamma^*$.
Now, $\gamma^*$
operates non-trivially on $\pi^*(E_t)$. This is because $\gamma$ is
not an $E_t$-automorphism of $C_t$. Because $\gamma^*$ is an involution,
it operates as $-\id$ on
$\pi^*(E_t)$. Thus $\pi^*(E_t)$ lies in the kernel of $1 +
\gamma^*$, i.e.\ it lies in the kernel of $c_* : J_{C_t} \longrightarrow
J_{D_t}$. This implies that $c^*(J_{D_t})$ lies in $\ker(\pi_*) = A_t$, the
complement of $\pi^*(E_t)$ under the Rosati involution.

Let $\tau := c_* \tau^*_{\ell} c^* \in \End^0_{k[\zeta_\ell]}((J_{D_t})_{k[\zeta_\ell]})$. We are interested in the minimal polynomial of $\tau$ and the Galois action on $\mathbb{Q}[\tau]$.
\footnote{If $A$ is some abelian variety over some field $K$, $\ell$ a prime $\neq \char(K)$, and $\alpha$ is some endomorphism on $A$, then the minimal polynomial of $\alpha$ in its operation on $V_\ell(A)$ lies in $\mathbb{Z}[T]$, in particular, it is equal to the minimal polynomial of $\alpha$ in the $\mathbb{Q}$-algebra $\mathbb{Q}[\alpha]$. We refer to this polynomial as the \emph{minimal polynomial} of $\alpha$.

This follows by induction on the degree of the minimal polynomial $m_\alpha$ of $\alpha$ in its operation on $V_\ell(A)$. Indeed, let $h$ be the product of all irreducible divisors of $\chi_\alpha$, the characteristic polynomial of $\alpha$. As $\chi_\alpha \in \mathbb{Z}[T]$ (see \cite[\S 19, Theorem 4]{Mu}), $h$ has the same property. Now, $h|m_\alpha$, and the minimal polynomial of $h(\alpha)$ in its operation on $V_\ell(A)$ is $\frac{m_\alpha}{h}$ which lies in $\mathbb{Z}[T]$ by induction assumption.}

The homomorphism $c^*$ induces an isogeny between $J_{D_t}$ and $c^*(J_{D_t}) = (\id + \gamma^*)(J_{C_t})$. In fact, $c_* c^* = 2 \, \id$ and $c^* c_*|_{c^*(J_{D_t})} = 2 \, \id$. This implies that we have an isomorphism of rings (with unity) and Galois modules
\[
\End^0_{k[\zeta_\ell]}((J_{D_t})_{k[\zeta_\ell]}) \longrightarrow \End^0_{k[\zeta_\ell]}(c^*(J_{D_t})_{k[\zeta_\ell]}), \; \;\; \alpha \mapsto \frac{1}{2} c^* \alpha c_*|_{(J_{D_t})_{k[\zeta_\ell]}}. \]
\pagebreak
Under this isomorphism, $\tau$ corresponds to $\frac{1}{2}c^* \tau c_*|_{c^*(J_{D_t})} = $\linebreak 
$\frac{1}{2} c^* c_*  \tau_\ell^* c^* c_*|_{c^*(J_{D_t})} =  \frac{1}{2}(\id + \gamma^*) \, \tau_\ell^* \, (\id + \gamma^*)|_{c^*(J_{D_t})} = (\tau_\ell^* + {\tau_\ell^*}^{-1}) \, \frac{1}{2} (\id + \gamma^*)|_{c^*(J_{D_t})} = (\tau_\ell^* + {\tau_\ell^*}^{-1})|_{c^*(J_{D_t})}$.
\footnote{In particular, $\tau_\ell^* + {\tau_\ell^*}^{-1}$ restricts to an endomorphism of $c^*(J_{D_t})$. This also follows from the fact that $\tau_\ell^* + {\tau^{*}_\ell}^{-1}$ and $\id +\gamma^*$ commute. The calculations in \cite[3.1]{ttv} are not necessary to prove this.}

Now, $\mathbb{Q}[\tau^*_\ell] \leq \End^0_{k[\zeta_\ell]}((A_t)_{k[\zeta_\ell]})$ is isomorphic to $\mathbb{Q}[\zeta_\ell]$ with $\tau \longleftrightarrow \zeta_\ell$. This implies that the minimal polynomial of $(\tau^*_\ell + {\tau^*_\ell}^{-1})|_{A_t}$ is equal to the minimal polynomial of $\zeta_\ell + \zeta_\ell^{-1}$.
It follows that the minimal polynomial of $\tau$, i.e.\ the minimal polynomial of $(\tau_\ell^* + {\tau_\ell^*}^{-1})|_{c^*(J_{D_t})}$, is also equal to the minimal polynomial of $\tau_\ell + \tau_\ell^{-1}$. 
We conclude that $\mathbb{Q}[\tau]$ is isomorphic to $\mathbb{Q}[\zeta_\ell + \zeta_\ell^{-1}]$ with $\tau \longleftrightarrow \zeta_\ell + \zeta_\ell^{-1}$.

Let $k = \mathbb{F}_{\!q}$ for some power $q$ of $p$. Then under the above isomorphism $\mathbb{Q}[\tau] \simeq \mathbb{Q}[\zeta_\ell + \zeta_\ell^{-1}]$, the operation of the Frobenius on $\tau$ corresponds to $\zeta_\ell + \zeta_\ell^{-1} \mapsto \zeta_\ell^q + \zeta_\ell^{-q}$. Thus $\tau$ is defined over $K:=\mathbb{F}_{\!q}[\zeta_\ell + \zeta_\ell^{-1}]$ and over no subfield $\mu$ of $K|k$ with $\mu \subsetneq K$. Note that $\Gal(\mathbb{F}_{\!q}[\zeta_\ell + \zeta_\ell^{-1}]|\mathbb{F}_{\!q}) \simeq \langle q \rangle \leq (\mathbb{Z}/\ell\mathbb{Z})^*/\langle -1 \rangle$.

Let $i > 1$ such that $p^i$ is a generator of
$(\mathbb{Z}/\ell\mathbb{Z})^*/\langle -1 \rangle$. As stated above,
there exists some $t \in \mathbb{F}_{\!p^i} \backslash \{ \pm 2 \}$ such that $D_t$ is ordinary.
We can apply Proposition \nolinebreak \ref{Weil-res} to $J_{D_t}$,
$k=\mathbb{F}_{\!p^i}$, $K=\mathbb{F}_{\!p^{i(\ell-1)/2}}$ and $\tau$.

We obtain that
$J_{D_t}$ is isogenous to the Weil restriction (with respect to $K|k$)
of one ordinary elliptic curve over $K$. We have proven the following
proposition which is slightly stronger than Theorem \ref{theorem-char-neq-2}.

\begin{proposition}
Let $p$ and $\ell$ be odd prime numbers, let $i \in \mathbb{N}, i >1$,
such that $p^i$ is a generator of
$(\mathbb{Z}/\ell\mathbb{Z})^*/\langle -1 \rangle$. Then there exists
a hyperelliptic curve over $\mathbb{F}_{\!p^i}$ of genus
$\frac{\ell-1}{2}$ whose Jacobian variety becomes over $\mathbb{F}_{\!p^{i(\ell-1)/2}}$ isogenous to the power of one elliptic curve. In
fact, there exists such a curve over $\mathbb{F}_{\!p^i}$ whose Jacobian is isogenous the the Weil restriction with respect to
$\mathbb{F}_{\!p^{i(\ell-1)/2}}|\mathbb{F}_{\!p^i}$ of one ordinary elliptic curve.

\end{proposition}
\begin{remark}
In \cite{ttv}, the curves $D_t$ have already been studied in
characteristic \nolinebreak 0. There it is shown that for $\ell \neq
5$ the Jacobian of the generic curve $D^\ell_t$ over $\mathbb{Q}(t)$
is absolutely simple (see \cite[Corollary 6]{ttv}). 
We think that the same is true for the generic curve $D_t^\ell$ over
$\mathbb{F}_{\!p}(t)$ for any $p$. Note however that by our above
results, if $p$ is a generator of
$(\mathbb{Z}/\ell\mathbb{Z})^*/\langle -1 \rangle$, for infinitely many
$t \in \overline{\mathbb{F}}_p$, $J_{D^\ell_t}$ is completely
decomposable. 
\end{remark}

\begin{remark}
As above, let $p^i$ be a generator of $(\mathbb{Z}/\ell\mathbb{Z})^*/\langle -1 \rangle$, and let $t \in
\mathbb{F}_{\!p^i} \backslash \{ \pm 2 \}$. We have used the endomorphism $\tau$ on
$(J_{D_t})_{\mathbb{F}_{\!q}[\zeta_\ell + \zeta_\ell^{-1}]}$ to derive that
$\chi_{J_{D_t}} = f(T^{(\ell-1)/2})$ for some polynomial $f \in
\mathbb{Z}[T]$ of degree 2.
But this can also be proven in an alternative way. Note that $p^i$ being a generator of $(\mathbb{Z}/\ell\mathbb{Z})^*/\langle -1 \rangle$ is equivalent to $p^{2i}$ generating ${(\mathbb{Z}/\ell\mathbb{Z})^*}^{2}$.
It is well known
that the Dickson polynomial $D_\ell(x,a)$ is a permutation polynomial for
$\F_{\!q}$ if $\gcd(q^2-1,l)=1$ (see e.g.\ \cite{dickson}). 
So $D_\ell(x,a)$ is a permutation polynomial for all $\F_{\!p^{ij}}$ with
$\frac{\ell-1}{2} \nmid j$, and consequently, 
$\#D_t(\F_{\!p^{ij}})=p^{ij}+1$ for those j.
It follows that the $L$-polynomial
$$L(D_t,T)=\exp\left(\sum_{j=1}^\infty
(p^{ij}+1-\#D_t(\F_{\!p^{ij}})) T^j\right)$$
is a polynomial in $T^{(\ell-1)/2}$. Since $\chi_{J_{D_t}}$ is the
reciprocal polynomial of $L(D_t,T)$, the same holds for $\chi_{J_{D_t}}$.
\end{remark}

\begin{remark}
Instead of the curves $C_t$ and $D_t$, it might be
possible to use other
ordinary hyperelliptic curves whose Jacobians have suitable endomorphisms. 
For example, in \cite{mestre}, some algebraic families of
hyperelliptic curves with real multiplication are given. If one could
show that these families of curves are generically ordinary, they
would give rise to new examples of hyperelliptic curves whose
Jacobians are isogenous a power of some ordinary elliptic curve, and
thus to new examples of high rank elliptic curves. We have  
observed that indeed, some of these curves are ordinary.
\end{remark}

\bibliography{highrank-literatur}

\bibliographystyle{plain}
\vspace{1 cm}
\noindent
Irene I.\ Bouw, Claus Diem: Institut f\"ur Experimentelle Mathematik, Ellernstr. 29, 45326 Essen, Germany. \{bouw, diem\}@exp-math.uni-essen.de.\\[0.5 ex]
Jasper Scholten: ESAT/COSIC, K.U. Leuven, Kasteelpark Arenberg 10, 3001 Leuven-Heverlee, Belgium. jasper.scholten@esat.kuleuven.ac.be.
\end{document}